\documentstyle[12pt]{article} 

\font\tenmsb=msbm10
\font\sevenmsb=msbm7
\font\fivemsb=msbm5

\newfam\msbfam
\textfont\msbfam=\tenmsb
\scriptfont\msbfam=\sevenmsb
\scriptscriptfont\msbfam=\fivemsb
\def\Bbb#1{{\fam\msbfam #1}}

\font\teneufm=eufm10
\font\seveneufm=eufm7
\font\fiveeufm=eufm5
\newfam\eufmfam
\textfont\eufmfam=\teneufm
\scriptfont\eufmfam=\seveneufm
\scriptscriptfont\eufmfam=\fiveeufm

\makeatletter
\ifnum\@ptsize=0  \fi \ifnum\@ptsize=2
 \fi 

\newcommand\qed{{\hspace*{\fill}Q.E.D.\vskip12pt plus 1pt}}

\newcommand\sE{{\cal E}}
\newcommand\sF{{\cal F}}

\newcommand\sJ{{\cal J}}
\newcommand\sI{{\cal I}}

\newcommand\sL{{\cal L}}

\newcommand\sN{{\cal N}}
\newcommand\sO{{\cal O}}

\newcommand\sT{{\cal T}}

\newcommand\grG{\Gamma}

\newcommand\grk{\kappa}

\newcommand\hotimes{\hat \otimes}

\newcommand\rat{{\Bbb Q}}
\newcommand\real{{\Bbb R}}

\newcommand\normB[2]{{{\cal N}_{{#1}\vert{#2}}}}
\newcommand\conormB[2]{{{\cal N}^*_{{#1}\vert{#2}}}}

\newcommand\sym[1]{{\rm \widehat S}^{{#1}}} \newcommand\symm[1]{{\rm S}^{{#1}}}

\newcommand\oline[1]{{\overline {#1}}}

\def\rank{{\mathop{\rm rank\,}\nolimits}} \def\cod{{\mathop{\rm cod}\nolimits}}

\def\sing{{\mathop{\rm Sing\,}\nolimits}}

\newcommand\proof{{\noindent\bf Proof.\ }}

\newtheorem{theorem}{Theorem}[section]
\newtheorem{lemma}[theorem]{Lemma}
\newtheorem{corollary}[theorem]{Corollary}

\newtheorem{proposition}[theorem]{Proposition}

\newtheorem{re}[theorem]{Remark} \newtheorem{pargrph}[theorem]{}
\newtheorem{examp}[theorem]{Example}

\newenvironment{remark}{\begin{re}\em}{\end{re}}

\newenvironment{prgrph*}[1]{\indent\begin{pargrph}{\bf #1.}\em\
}{\end{pargrph}}

\begin{document}
\title{Kodaira Dimension of Subvarieties} \author{Thomas Peternell\and
Michael Schneider\and Andrew J. Sommese} \date{March 26, 1998} \maketitle

\section*{Introduction} In this article we study how the birational
geometry of a normal projective variety $X$ is influenced by a normal
subvariety $A \subset X.$ One of the most basic examples in this
context is provided by the following situation. Let $f:X\to Y$ be a
surjective holomorphic map with connected fibers between compact connected
complex manifolds. It is well known (see, e.g., \cite{Ueno}) that given a
general fiber $A$ of $f$ we have $$
\grk(X)\le \grk(A)+\dim Y.
$$
This article grew out of the realization that this result should be true with $\dim Y$ replaced by the codimension $\cod_X A$ for a
pair $(X,A)$ consisting of a normal subvariety $A$ of a compact normal
variety $X$ under weak semipositivity conditions on the normal sheaf of $A$
and the weak singularity
condition $\cod_A (A\cap\sing X)\ge 2$.
We shall now state our main results in the special case of a submanifold
$A$ in a projective
manifold $X$ and we also simplify the semipositivity notion.

\begin{theorem} Let $X$ be a projective manifold and $A$ a compact
submanifold. Then
$$ \grk (X) \leq \grk (A) + \cod_X A$$
if one of the following conditions is satisfied \begin{enumerate}
\item some symmetric power $S^m\sN_{A \vert X}$ of the normal bundle has
global sections which generate the
bundle almost everywhere;
\item $\sN_{A \vert X}$ is nef and $A$ has a good minimal model via
contractions and flips;
\item $\sN_{A \vert X}$ is nef and $\dim A \leq 3;$ \item $\cod_X A = 1$,
$\sN_{A \vert X}$ is in the closure of the effective cone of $A$
and $A$ has a good minimal model via contractions and flips. \end{enumerate}
\end{theorem}

Recall that a normal projective variety $Z$ is a good minimal model if $Z$
is $\Bbb Q-$factorial with at most terminal singularities
and some multiple $mK_Z$ is generated by global sections.

The above results actually hold more generally, e.g., with $A$ a normal
projective subvariety of a normal projective variety $X$ such that $\cod_A
(A\cap\sing X)\ge 2$. In this case it is necessary to use the arithmetic
Kodaira dimension of $X$ and $A$. Moreover the effectivity and nefness
assumptions can be weakened in the following way; we assume again $X$ and
$A$ smooth: some symmetric power $S^m\sN_{A \vert X}$ has a decomposition
$S^m\sN_{A \vert X} \simeq A \otimes B$ with $A$ a generically spanned and
$B$ a ``generically nef" vector bundle. Generically nef means the
following: there is a
Zariski open set $U \subset X$ with $\cod_X (X \setminus U) \geq 2,$ such
that $B \vert U$ is nef, i.e., $B \vert C$ is nef for every compact curve
$C \subset B.$

For technical reasons we have formulated most part of the paper in a
non-compact setting,
namely for normal pairs $(X,A).$ Here $X$ is a normal variety and $A
\subset X$ a normal
subvariety such that $A \cap \sing (X)$ has codimension at least 2 in $A.$
We require the
existence of normal compactifications $\oline X$ and $\oline A$ (not necessarily
$\subset \oline X$) such that the boundary components have codimension at
least 2. The
reason for using this category is that at some point we have to perform
surgery in
codimension 2 (flips) and then this language seems appropriate.

This paper being almost completed, Michael Schneider died in a tragic
accident. The scientific community lost a very active
mathematician; we lost also a very good friend. We dedicate this work to
his memory.

\section{Preliminaries}For most purposes of this article singular varieties
are not much
harder to deal with than smooth varieties, except that some care must be
taken with definitions.
A normal variety is a connected normal quasi-projective scheme over the
complex numbers (or a
connected normal quasi-projective complex space).

Let
$\sF$ be a coherent sheaf on a complex algebraic variety $X$. Coherent
sheaves are of course
always understood to be algebraic.
By $\sym t \sF$ we denote $(\symm t\sF)^{**}$, the double dual of the
$t$-th symmetric power, $\symm t\sF$, of $\sF$.
When $\sL$ is of rank one reflexive sheaf on a normal variety $X$, we
often, for $t\ge 1$, denote $\sym t \sL$ by
$t\sL$.

If $A$ and $B$ are reflexive sheaves, we define $$ A \hat \otimes B = (A
\otimes B)^{**}.$$

The following is left to the reader.
\begin{lemma}Let $\sF$ denote a reflexive sheaf on a normal variety $V$.
Assume that there is
an embedding $i:V\to \oline V$ of $V$ as a Zariski open set in a compact
normal variety $\oline V$. Then $\sF$ extends to a reflexive sheaf $\oline
\sF$ on $\oline V$. If $\cod_{\oline V}\oline V-V\ge 2$ then $\sym t\sF$
extends to $\sym t\oline \sF$ for all t$\ge 0$, and moreover $h^0(\sym
t\sF)=h^0(\sym t\oline\sF)$.
Furthermore $i_*\sym t\sF\cong \sym t\oline\sF$. \end{lemma}

Given a rank one reflexive sheaf $\sL$ on a normal variety $V$ the {\em
Kodaira dimension, $\grk (\sL)$ of} $\sL$, is defined as follows
\begin{enumerate}
\item $\kappa (\sL) = -\infty$, if $h^0(\sym t\sL)=0$ for all $t\ge 1$;
\item $\kappa (\sL) = 0$, if $h^0(\sym {t_0}\sL)= 1$ for some $t_0\ge 1$
and $h^0(\sym t\sL)\le 1$ for all $t\ge 1$; \item $\kappa (\sL)$ is a
positive integer $k$, if $\displaystyle\oline \lim_{\{t\ge 1\}}
\frac{h^0(\sym t \sL)}{t^k}$ is a positive real number, where $\oline \lim$
denotes the $\mbox{\rm limsup}$; \item $\kappa (\sL) = \infty$ otherwise.
\end{enumerate}

We have the following corollary of the above lemma. \begin{lemma}Let $\sL$
denote a reflexive rank one sheaf on a normal variety $V$. Assume that
there is an embedding $i:V\to \oline V$ of $V$ as a Zariski open set in a
compact normal
variety $\oline V$. Let $\sL$ extend to the reflexive sheaf $\oline \sL$ on
$\oline V$. If $\cod_{\oline V}(\oline V-V)\ge 2$ then for any $t\ge 1$
$$\grk(\sL)=\grk(t\sL)=\grk(t\oline\sL)=\grk(\oline\sL)<\infty$$
\end{lemma}

By a {\em normal pair} $(X,A)$ we mean a pair $(X,A)$ where:
\begin{enumerate} \item $A$ is a normal subvariety of a normal variety $X$;
and
\item there is an embedding $A\subset \oline A$ of $A$ as a Zariski open
set in a normal
compact projective variety $\oline A$ in such a way that $\cod_{\oline
A}(\oline A-A)\ge 2$; and \item there is an embedding $X\subset \oline X$
of $X$ as a Zariski open set in a normal
compact projective variety $\oline X$ in such a way that $\cod_{\oline
X}(\oline X-X)\ge 2$; and \item $\cod_A (A\cap\sing X)\ge 2$.
\end{enumerate}
Note we do not require that $\oline A$ equals the closure of $A$ in $\oline
X$. The main
example will be given by a normal projective variety $X$ and a normal
projective subvariety $A$ such that $\cod_A(A \cap \sing X) \ge 2.$

The condition that $\oline X$ is projective is needed for the basic
inequality Theorem \ref{basicInequality}. The condition that $\oline A$ is
projective is needed for Theorem \ref{realEffectivityLemma}.
The codimension two conditions are needed for finiteness results. Indeed it
is very easy to see
that without some such conditions
the results are at best meaningless, e.g., let $A$ be a smooth compact
curve and let
$X := A\times (C-\{x\})$ for some point $x$ on a smooth compact curve $C$.

Some further notations:
\begin{enumerate}
\item Given a normal pair $(X,A)$, we say that {\em a reflexive sheaf $\sF$
on $X$ is adapted to $(X,A)$},
if $A \not \subset \sing \sF.$
\item A coherent sheaf $\sF$ is nef on a normal variety $V$ if $\sF$ is nef
on every compact
irreducible curve $C \subset V,$ i.e. $\sO_{\Bbb P (\sF \vert C)}(1) $ is
nef on $\Bbb P (\sF \vert C).$
\item For a coherent sheaf $\sF$ we say that $\sF$ has rank $r$ if there is
a Zariski open and dense set on which it is a locally free sheaf of rank
$r$. In this case
we define we define the determinant of $\sF$ by  ${\rm det}(\sF) =
(\bigwedge ^r(\sF))^{**}.$
\end{enumerate}

\section{The basic inequality}\label{basicInequalitySection}The following
is the basic inequality underlying this paper.
\begin{theorem}\label{basicInequality}Let $(X,A)$ denote a normal pair and
let $\sL$ be
a reflexive rank
one coherent sheaf
on $X$ adapted to $(X,A)$. Then there is a positive integer $c$ such that
for all $t\ge 0$
$$
h^0(t\sL)\le \sum_{k=0}^{ct} h^0(\sym k\conormB A X \hotimes t\sL_A). $$ In
particular, if
we have
$$h^0\left(\sym k\conormB A X \hotimes t\sL_A\right)\le C\left(\rank \symm
k\conormB A X \right) t^{a}$$
for some positive
constants $C, a$ that do not depend on $t, k$, then we have $$\grk(\sL)\le
\cod_XA+a.$$
\end{theorem}
\proof Let $\sJ_k = \sI_A^k.$
The essential point is to show that there is a positive integer $c$ such
that $h^0(t\sL\otimes \sJ_k)=0$ for $k>ct$. Since there are projective
varieties $\oline A,\oline X$ such that $\cod_{\oline X}(\oline X-X)\ge 2$
and $\cod_{\oline A}(\oline A-A)\ge 2$, it follows that we can assume
without
loss of generality that
$A$ and
$X$ are compact.

If $A$ is a divisor then choose a very ample curve $C\subset X$, i.e., the
intersection of $\dim X-1$ very ample divisors on $X$. By choosing the very
ample divisors generically we can assume that \begin{enumerate}
\item $C$ is a smooth connected curve lying in $X-\sing X$ and meeting $A$
transversely in points lying in $A-\sing A$; and
\item the restriction
$\sL_C$ is invertible.
\end{enumerate}
Since
$$\deg(t\sL\otimes \sJ_k)_C =t\deg\sL_C-kA\cdot C,$$ we see that
$h^0((t\sL\otimes \sJ_k)_C)=0$ for $\displaystyle
k>\left(\frac{\deg\sL_C}{A\cdot C}\right)t$. Since $C$ is very ample this
implies that
$h^0(t\sL \otimes \sJ_k)=0$ for
$\displaystyle k> \left(\frac{\deg \sL_C}{A\cdot C}\right)t$. So let $c :=
\displaystyle k>\left(\frac{\deg\sL_C}{A\cdot C}\right)t$ and our
inequality
follows by power series expansion.

If $A$ is not a divisor we proceed as follows. Again it is sufficient to
show $$ H^0(C,\sym l\conormB A X \otimes t\sL_C) = 0 \eqno (*)$$ for $k >
ct,$ with $c$ not depending on the individual curve $C.$ To verify (*) we
choose $C$ again general, blow up $X$ in a neighborhood of $C$ where $X$
and $A$ are smooth and apply the old argument. \qed

\section{$\rat$-effective and $\real$-effective sheaves}We need the rank
$>1$ version of a
$\rat$-effective divisor, i.e., of a $\rat$-effective rank one coherent
sheaf. A coherent sheaf $\sF$ on a normal variety $V$ is said to be {\em
generically
spanned} if the global sections $\grG(\sF)$ span $\sF$ over a dense Zariski
open set of $V$.
We say that a coherent sheaf $\sF$ on a normal variety $V$ is {\em
$\rat$-effective} if there is a positive integer $N>0$ such that $\sym
N\sF$ is generically spanned. Note that for a line bundle, $\rat$-effective
agrees with the usual notion.

The proofs of the following lemmas are immediate.

\begin{lemma}\label{QeffectiveRestriction}Let $\sF$ be a reflexive sheaf on
a normal variety $V$.
Then the following are equivalent:
\begin{enumerate}
\item $\sF$ is $\rat$-effective;
\item $\sF_U$ is $\rat$-effective for every dense Zariski open set
$U\subset V$ with $\cod_V(V-U)\ge 2$;
\item $\sF_U$ is $\rat$-effective for some dense Zariski open set $U\subset
V$ with $\cod_V(V-U)\ge 2$.
\end{enumerate}
\end{lemma}

\begin{lemma}Let $\sF$ be a $\rat$-effective coherent sheaf on a normal
variety $V$. Then any
coherent quotient sheaf $\sT$ of $\sF$ on $V$ is also $\rat$-effective.
Further $\sym b\sF$ is
$\rat$-effective for all $b>0$.
\end{lemma}

We need a generalization of the tensor product of a $\rat$-effective with a
nef line bundle.
We say that a coherent sheaf $\sF$ on a normal variety $V$ is {\em
$\real$-effective} if there is a positive integer $N>0$ such that
$\displaystyle \sym N\sF \cong A \hotimes B = (A\otimes B)^{**}$ where
\begin{enumerate}
\item $A$ is a generically spanned reflexive sheaf on $V$; and \item $B$ is
a reflexive sheaf on $V$ with $B_U$ a nef locally free sheaf for some dense
Zariski open set $U\subset V$ with $\cod_V(V-U)\ge 2$. \end{enumerate}

\bigskip \noindent $\real$-effectivity should be a notion generalizing the
property of a
line bundle to be in the closure of the effective cone. However it is not
clear whether for
line bundles the two notions coincide. Line bundles which are in the
closure of the effective
cone are studied in \cite{DPS}, where they are called pseudo-effective.

\bigskip \noindent The proof of the following lemma is again immediate.

\begin{lemma}\label{realEffectiveEquivalence}Let $\sF$ be a torsion free
coherent sheaf on a
normal variety $V$. The following are equivalent: \begin{enumerate}
\item $\sF$ is {\em $\real$-effective};
\item There is a positive integer
$N>0$, a dense Zariski open set $U\subset V$ with $\cod_V(V-U)\ge 2$, and
reflexive
sheaves $A,B$ on $V$ such that \begin{enumerate} \item $U$ is smooth;
\item $A_U$, $B_U$, $\sF_U$ are locally free; \item $\symm N\sF_U\cong
A_U\otimes B_U$;
\item $A_U$ is a generically spanned locally free sheaf on $V$; and \item
$B_U$ is a nef locally free sheaf. \end{enumerate} \end{enumerate}
\end{lemma}

\begin{theorem}\label{realEffectivityLemma}Let $\sF$ be a torsion free
coherent sheaf on a
normal variety $V$ and assume that $\sF^*$ is $\real$-effective. Assume
that $V$ is a Zariski
open subset of a projective variety $\oline V$ such that $\cod_{\oline
V}(\oline V-V)\ge 2$,
e.g., $V$ is projective.
If $s\in \grG(\sF)$ is not identically zero, then $s$ is non-vanishing on
any Zariski open set
$U\subset V$ satisfying the conditions of Lemma
\ref{realEffectiveEquivalence}. In particular, $h^0(\sF)\le \rank \sF$.
\end{theorem} \proof By replacing $V$ by the open set $U$ satisfying the
conditions of Lemma \ref{realEffectiveEquivalence} for $\sF^*$, we can
assume without loss of generality that $V$ is smooth, $\sF$ is locally
free; and that there is a positive integer
$N>0$ and locally free sheaves $A,B$ on $V$ such that \begin{enumerate}
\item $\symm N\sF^*\cong A\otimes B$;
\item $A$ is a generically spanned locally free sheaf; and \item $B$ is a
nef locally free sheaf;
\item $V$ is a Zariski open subset of a projective variety $\oline V$ such
that $\cod_{\oline V}(\oline V-V)\ge 2$. \end{enumerate} Replacing $\sF$
with $\symm N\sF$ and $s$ with $\symm Ns $ we can assume without loss of
generality that
$N=1$.

We must show that $s$ does not vanish on $V$.

By assumption $s$ gives rise to a non-zero map $\sF^*\to \sO_V$. Assume
that there was a
point $x\in V$ with $s(x)=0$. Since $V$ is a Zariski open subset of a
projective variety
$\oline V$ such that $\cod_{\oline V}(\oline V-V)\ge 2$, we can choose a
smooth irreducible
projective curve $C\subset V$ with $x\in C$, with $s_C$ a not identically
vanishing section
of $\sF_C$ that vanishes at $x$. If $\sF^*\to \sO_C\to 0$, then we are
done. If this is not true
then we have $\sF^*\to L^*\to 0$ for an ample line bundle $L$ on $C$. On
the other hand,
$$ \sF^*_C \cong A_C \otimes B_C.$$
If $C$ is general (through $x$), then $A_C$ is generically spanned, hence
nef, and $B_C$ is
nef, too. Since a quotient of a nef bundle has to be nef, we obtain a
contradiction.
\qed
The above result immediately gives a number of strong consequences.
\begin{theorem}Let $(X,A)$ be a normal pair. Assume that the normal sheaf
$\normB A X$ is $\real$-effective. Assume that $\sL$ is a rank one
reflexive sheaf on $X$ adapted to $(X,A)$. If $\kappa (\sL_A^{*}) > 0$ then
$\kappa (\sL)= -\infty$.
\end{theorem}
\proof By Theorem \ref{basicInequality} it suffices to show that $$h^0(A,
t\sL_A \hotimes \sym k\conormB A X)=0$$ for $t>0$ and $ k<c t$. Assume that
we have a nonzero section
$$\sigma \in H^0(A, t\sL_A \hotimes \sym k\conormB A X).$$ By assumption,
for some $m>0$ we have a non trivial section $\tau\in H^0(A, {-m}\sL_A)$
with zeros. This leads to an inclusion ${mt}\sL_A\to \sO_A$ with zeros and
altogether we obtain a nontrivial section with a divisor of zeros in
$H^0(A, \sym {km}\conormB A X)$. This contradicts Theorem
\ref{realEffectivityLemma}.
\qed
\begin{corollary}
Let $X$ be a normal projective variety and let $A\subset X$ be a smooth
projective
curve with $\real$-effective normal bundle, e.g., the normal sheaf $\normB
A X$ is nef or
generically spanned. Assume $A \cap \sing X = \emptyset.$ If $\kappa (X) =
\kappa (K_X) \geq 0$
then $K_X\cdot A \geq 0$.
\end{corollary}
Another variant of this is the following result. \begin{theorem} Let $X$ be
a normal projective variety and let $A\subset X$ be a smooth projective
curve with
$\real$-effective normal bundle, e.g., the normal sheaf $\normB A X$ is nef
or generically
spanned. Assume $A \cap \sing X = \emptyset.$ If $X$ is of general type
then $K_X\cdot A> 0$. \end{theorem} \proof By the previous corollary we
only need to exclude the case $K_X\cdot A= 0$. Assume that $K_X\cdot A= 0$.
Note that given a torsion free sheaf on $A$ whose dual is
$\real$-effective,
we have
$$h^0(A, \sF)\leq \rank \sF$$
by the last conclusion of Theorem \ref{realEffectivityLemma}. This will be
applied to
$(tK_X \vert A)\otimes \sym k\conormB A X$ in the standard estimate
$$h^0(X, tK_X)\leq\sum_{k=0}^{c t}h^0((tK_X \vert A)\otimes \sym k\conormB
A X)\leq
\sum_{k=0}^{c t}\rank(\sym k\conormB A X).$$ For large $t$ this grows like
$t^{n-1}$, $n= \dim X$. Hence $\kappa (X)\leq n-1$, contradicting
our hypothesis. \qed

\bigskip \noindent We now prove three lemmas which will be important for
the proof of the
main results in the next section.

\begin{lemma}Let $\sL$ be a reflexive $\rat$-Cartier rank one coherent
sheaf on a normal
variety $V$.
Assume that $V$ is a Zariski open subset of a projective variety $\oline V$
such that
$\cod_{\oline V}(\oline V-V)\ge 2$. If $\sL$ is semiample and $\sF^*$ is a
$\real$-effective
sheaf on $V$, then we have for $t$ sufficiently divisible $$h^0\left(\sF
\hotimes t\sL\right)\le C\left(\rank\sF\right)t^{\grk(\sL)} $$ where
$C$ is a positive constant that depends only on $(V,\sL)$ and a data of
section and neither
on $\sF$ nor on $t.$ The data of sections is described as follows. Take the
smallest number
$t_0$ such that $t_0\sL$ is locally free and spanned. Let $f: V \to W$ be
the associated morphism
and write $\sL = f^*(\sL').$ Now choose $s'_1, \ldots, s'_d, d = \dim W,$
general such that the
common intersection is finite. Then $s_i = f^*(s'_i)$ form the data of
sections we fix. \end{lemma}

\proof By passing to a suitable multiple of $\sL$ if necessary, we may
assume that $\sL$ is
locally free, that $\sL$ is already spanned and moreover that $V$ is
compact. We choose
$s$ minimal.

Let $f: V \to W$ be the Stein factorization of the morphism associated to
the linear system
$\vert \sL \vert.$ Then $\sL = f^*(\sL');$ we may assume that $\sL'$ is
very ample.
We shall proceed by induction on $d = \dim W.$

The case $d = 0$ is obvious because then $\sL = \sO_V.$

So suppose $\dim W = d > 0.$ We fix a smooth member $H \in \vert 2L'
\vert.$ Let $V_t =
f^{-1}(tH).$ Then we have an exact sequence for $t > 0$ $$ 0 \to
H^0(V,\sL^t \otimes \sI_{V_t} \otimes \sF) \to H^0(V,\sL^t \otimes \sF) \to
H^0(V_t,\sL^t \otimes \sF \vert V_t).$$
Now our claim will be a consequence of the two following assertions
\begin{enumerate} \item $H^0(V,\sL^t \otimes \sI_{V_t} \otimes \sF) = 0$
and
\item $h^0(V_t,\sL^t \otimes \sF \vert V_t) \leq C (\rank \sF) \ t^{\grk
(\sL)}.$
\end{enumerate}

For the proof of (1) notice first
$$ \sL^t \otimes \sI_{V_t} \otimes \sF = f^*(\sL^{'(-t)}) \otimes \sF.
\eqno (*)$$
Now suppose that $s$ a non zero section of $\sL^t \otimes \sI_{V_t} \otimes
\sF.$
Let $U \subset V$ be a big open set in the sense of (3.3). Let $C \subset
U$ be a general compact irreducible curve.
Take a non-zero section $s' \in H^0(C,f^*(\sL^{'t} \vert f(C))).$ Since
$\sL'$ is
ample, $s'$ has zeros. Then, using (*), $s \otimes s' \in H^0(C,\sF \vert
C))$ is
a section with zeros contradicting (3.4).

As to the proof of (2) we use the fact that $V_t$ is the $t-$th
infinitesimal neighborhood
of $V_1 = V_H = f^{-1}(H).$ Hence
$$ h^0(V_{tH},\sL^t \otimes \sF \vert V_{tH}) \leq \sum_{\mu = 0}^{t-1}
h^0(V_H,\sL^{t-2\mu}
\otimes \sF).$$
By induction there is a constant independent of $t$ and $\sF$ such that $$
h^0(V_H,\sL^{t-2\mu} \otimes \sF) \leq C (\rank \sF) (t-2\mu)^{\grk (\sL
\vert {V_h})}$$
if $\dim W \geq 2$ and
$$ h^0(V_H,\sL^{t-2\mu} \otimes \sF) \leq C (\rank \sF) \ {\rm deg} \sL' $$
if $\dim W = 1.$ In this last case we substitute $C$ by $C {\rm deg} \sL'.$
Adding up and observing $\grk (\sL_H) = \grk (\sL) - 1,$ we obtain $$
h^0(V_{tH},\sL^t \otimes \sF \vert V_{tH}) \leq C (\rank \sF) \ t^{\grk
(\sL)}.$$
\qed

\bigskip \noindent Actually one can prove the last lemma for every
reflexive sheaf of rank
one but we do not need this. It seems however to be an interesting problem
whether 3.8 holds under
more general assumptions, e.g. when $\sL$ is nef.

\bigskip Next we investigate the behavior of $\Bbb R$-effective sheaves
under birational rational which are
isomorphisms out side sets of codimension ate least 2 (main examples are flips).
We fix the following situation. Let $V$ and $V'$ be normal varieties, $Y
\subset V$ and
$Y' \subset V'$ algebraic subsets of codimension at least 2 and let
$\varphi: V \rightharpoonup V'$
be birational such that $\varphi : V \setminus Y \to V' \setminus Y'$ is an isomorphism.
Let $\psi : W \to V$ be birational such that the induced map $\tau : W \to
V'$ is a morphism.
If $\sF$ is a torsion free sheaf on $V$, we define $$ \varphi_*(\sF) =
\tau_*(\psi^*(\sF))^{**}.$$ In other words, $\varphi_*(\sF)$ is the
"canonical" extension of $\varphi_*(\sF \vert V \setminus Y).$
In this context we can state

\begin{lemma} Assume that $\sF^*$ is $\real$-effective. Then
$(\varphi_*\sF)^*$ is $\real$-effective. \end{lemma}

\proof Write
$$ \sym m \sF = A \hotimes B$$
for a suitable positive integer $m,$ where $A$ is $\Bbb Q$-effective and
$B$ nef on $U \subset V$ with $\cod_V U \geq 2.$ Let $A' = \varphi_*(A)$
and $B' = \varphi_*(B).$ Then $B'$ is again
$\Bbb Q$-effective and since we can assume $U \cap Y = \emptyset, $ $A'$ is
again $\Bbb Q$-effective.
Since obviously
$$ \sym m \varphi_*(\sF) = A' \hotimes B',$$ the claim follows. \qed

\bigskip \noindent Concerning global sections we have in the same situation
with reflexive
sheaves $\sL$ and $\sF$ on $V:$

\begin{lemma}
$$ H^0(V,t\sL \hotimes \sF) = H^0(V',\varphi_*(t\sL) \hotimes \varphi_*(\sF)).$$

\end{lemma}

\proof This is clear since both sheaves coincide on $V \setminus Y = V'
\setminus Y' $ and
since $Y$ and $Y'$ have codimension at least 2. \qed

\section{Main Results}If $X$
is a normal variety we again let
$\grk(X)$ denote the arithmetic Kodaira dimension of $X$, i.e.,
$\grk(K_X)$. Note that $\grk(X)=\grk(U)$ for any Zariski open set $U\subset
X$ such that $\cod_X(X-U)\ge 2$.

\begin{theorem}Let $(X,A)$ be a normal pair. If $\normB A X$ is
$\rat$-effective, then
$$\grk (X) \le \grk(A) + \cod_X A.$$
\end{theorem}
\proof We are going to apply (2.1) with $\sL = K_X.$ First notice that
since $\sN_{A \vert X}$
is $\Bbb Q $-effective, we have an inclusion $\sym {k_0} \sN_{A \vert X}^*
\subset \sO_A^M$ for
a suitable positive integer $k_0.$
Hence
$$ h^0(\sym m \sN_A^* \hotimes t\sL_A) \leq \dim \sym {k_0} H^0(S^{k_0}
\sN_A^* \hotimes \sL_A)$$
$$ \leq h^0(\sym {k_0} \sym m \sN^*_A \otimes k_0t\sL_A) \leq C \ (\rank \sym m
\sN^*_A) \ h^0(k_0t\sL_A) $$ $$ \leq C' \ (\rank \sym m \sN^*_A) \ t
^{\kappa (\sL_A)}.$$
Thus we obtain
$$ \grk (X) \leq \grk (K_X \vert A) + \cod_X A.$$ Now $K_X \vert A = K_A
\otimes {\rm det}N_A^*$ at least outside a set of codimension
$\geq 2.$ Since ${\rm det}N_A$ is $\Bbb Q$-effective, we conclude
(similarly as above)
that $$\grk(K_X \vert A) \leq \grk (A)$$ and our claim follows.
\qed

\bigskip \noindent Actually most parts of (4.1) hold in more generality, we
have

\begin{theorem}Let $(X,A)$ be a normal pair. Let $\sL$ be a torsion free
sheaf of rank 1 on
$X$ and assume $\cod_A (A \cap \sing \sL) \geq 2.$ Assume that $\sN_{A
\vert X}$ is
$\Bbb Q$-effective. Then
$$ \grk (\sL) \leq \grk (\sL_A^{**}) + \cod_X A.$$ \end{theorem}

\bigskip \noindent We next discuss the case that $A$ has $\Bbb R$-effective
normal bundle,
e.g., nef normal bundle. To work out the difficulties with this case, let
us assume for
simplicity that $A$ is a Cartier divisor in $X.$ A main point in (4.1) was
that some
power of $N_A$ has a section and that therefore the dual is a subsheaf of
$\sO_A.$
If say $N_A$ is nef, then it might happen that no power of $N_A$ has a
section, so we
cannot argue in this way. To make this more concrete, suppose that $X$ is a
projective
manifold, that $\sL$ and $\sF$ are line bundles on $X$ and that $\sF^*$ is
nef. Intuitively one would say that $h^0(\sF \otimes L) \leq h^0(L)$
Of course this is false for trivial reasons. Take e.g., $X$ to be an abelian
variety, $\sL \in {\rm Pic}^o(X)$ non-torsion and $\sF = \sL^*.$ On the
other hand we need only
$h^0(\sF \otimes \sL^t) \leq C \rank h^0(\sL^t)$ asymptotically.

Lemma 3.8 shows that this is indeed true if $\sL$ is semi-ample, in
particular for $\sL = K_X$
when $X$ is a good minimal model.
The lemmas 3.9 and 3.10 show that if our original $A$ has a good minimal model
in the strong sense, i.e., via divisorial contraction of an extremal ray
(i.e. contractions of
a prime divisor) or flips (see e.g., \cite{KMM}), then we can go to $A'$ and argue there. Thus we have

\begin{theorem}Let $(X,A)$ be a normal pair. Assume $\oline A$ has at worst
terminal singularities
and that $\oline A$ admits a good minimal model in the strong sense. If
$\normB A X$ is
$\real$-effective, then $$\grk(X)\le \grk(A)+\cod_X A.$$ \end{theorem}
\proof
We may assume $(X,A)$ projective. Let $\varphi : A \rightharpoonup A'$ a
sequence of divisorial contractions and flips such that $A'$ is a good
minimal model.
By our basic inequality we have $$ h^0(tK_X) \leq \sum_{k=0}^{ct} h^0(\sym
k \sN^*_{A \vert X} \hotimes tK_A \hotimes {\rm det} t\sN^*_{A \vert X}).$$
Let
$$ \sF_{t,k} = \sym k\sN^*_{A \vert X} \hotimes t{\rm det} \sN_{A \vert
X}.$$ By (3.9) and (3.10) we conclude
$$ \sum_{k=0}^{ct} h^0(\sym k \sN^*_{A \vert X} \hotimes tK_A \hotimes {\rm
det}t\sN^*_{A \vert X}) = \sum_{k=0}^{ct} h^0(\varphi_*(\sF_{t,k}^*)
\otimes tK_{A'})) $$
$$\leq C \sum_{k=0}^{ct}
\rank \varphi_*(\sF^*_{t,k}) t^{\grk (A)}.$$ Here we have used the
semi-ampleness of $K_{A'}$ in order to be able to apply (3.8) and
moreover $\grk (A) = \grk (A').$ Note also that the constant $C$ is
independent of $t$ and $k$ by
(3.8). Since
$$ \rank \varphi_*(\sF^*_{t,k}) = \rank \sF_{t,k} \sim t^{\cod_X A - 1},$$
we conclude finally that
$$ h^0(tK_X) \leq C t^{\grk (A) + cod A}$$ asymptotically.
\qed

\bigskip \noindent

Parts of the proof actually show

\begin{theorem} Let $(X,A)$ be a normal pair and $\sL$ a reflexive sheaf on
$X.$ Assume that
$\sL$ is $\rat$-Cartier, that $\sL_A^{**}$ is semi-ample and adapted to
$(X,A).$ Assume that
$\sN_{A \vert X}$ is $\real$-effective. Then $$\grk (\sL)\leq \grk
(\sL_A^{**}) + \cod_X A.$$ \end{theorem}

\bigskip \noindent Without the semi-ampleness assumption on $\sL$ (4.3) is
however false if the normal bundle $\sN_{A \vert X}$ is, say, only nef. In
fact, let $C$ be an elliptic curve and let $\sF \in {\rm Pic}^0(C)$ be a
non-torsion point. Set $X = \Bbb P (\sO \oplus \sF)$ and let $A \subset X$
be the section
given by the epimorphism $\sO \oplus \sF \to \sO.$ Set $\sL = \sO_X(A).$
Then $\sL_A =
\sF^{-1}$ and therefore $\kappa (L_A) = - \infty.$ On the other hand
$\kappa (L) = 0.$
If however we define a refined Kodaira dimension $\tilde \grk$ substituting
$- \infty$ by $-1$ in the
definition, then this example does not work.

\begin{remark}It is worth noting that the same methods yield some easy but
very useful results in special situations. For example, assume that $A$ is
a positive dimensional connected compact complex submanifold of a connected
complex, but not necessarily compact, manifold $X$. Assume that
$h^0(S^k\sN^*_{A \vert X})=0$ for all $k>0$, e.g., assume that the normal
bundle of $A$ in $X$ is ample. Then given any holomorphic map $\phi:X\to Y$
from $X$ to a complex space $Y$ with $\dim\phi(A)=0$, it follows that
$\phi$ maps $X$ to a point. Thus if we have a globally generated
line bundle $L$ on a complex manifold $X$ and $L$ is trivial on a positive
dimensional compact submanifold with ample normal bundle, it follows that
$L$ is the trivial bundle.
\end{remark}

\bigskip \noindent {\bf Question} {\it \begin{enumerate} \item Assume
$\sN_{A \vert X}$ to be $\real$-effective and $\sL$ nef.
Is $$ \grk (\sL) \leq  \grk (\sL_A) + \cod_X A ?$$ \item Assume instead of
nefness of $\sL$ that $\sL$ is big and its canonical ring $\bigoplus _m
H^0(X,\sL^m)$
is finitely generated. Do we have
$$ \grk (\sL) \leq \grk (\sL_A) + \cod_X A ?$$ \end{enumerate} }

\bigskip \noindent It is easy to see that the answer to part (1) of the
question is yes for surfaces.

\bigskip \noindent By \cite{Mo,Mi,Ka} every threefold with only terminal
singularities
is either uniruled or has a good minimal model in the strong sense. Hence
we can state (we stick
to the most important case that $X$ and $A$ are projective).
\begin{corollary} \begin{enumerate} \item Let $(X,A)$ be a normal pair with
$X$ and  $A$ a projective
threefold such that
$\oline A$ has only terminal
singularities and $A$ is not uniruled. If $\sN_{A \vert X}$ is $\Bbb
R$-effective, then
$$ \grk (X) \leq \grk (A) + \cod_X A.$$
\item Let $A$ be a non-uniruled connected projective submanifold of a
projective manifold $X.$
If $\normB A X$ is nef and $\dim A\le 3$, then $\grk(X)\le \grk(A)+\cod_X
A$. \end{enumerate} \end{corollary}

\bigskip \noindent It remains to consider the case when $A$ is uniruled.
Here we consider only
the projective situation. It is necessary to make a slightly stronger
assumption on the
normal sheaf. On the other hand it is not necessary to assume anything on
the singularities
of $A.$

\begin{theorem} Let $(X,A)$ be a normal projective pair. Assume that $A$ is
uniruled and has
only terminal or canonical singularities. Assume furthermore that $\sym m
\sN_{A \vert X} \simeq A \hotimes B,$ where $A$ is a nef reflexive
sheaf and $B$ is a $\Bbb Q$-effective sheaf for some positive integer $m$.
Then $\grk (X) = - \infty.$ \end{theorem}

\proof
Let $(C_t)$ be a covering family of rational curves. Since $A \vert C_t$ is
nef (possibly
with torsion!) and $B \vert C_t$ is nef for general t and since moreover
$K_A \cdot C_t < 0$
(here we need the assumption on the singularities), it follows $$ H^0(\sym
k \sN^*_{A \vert X} \hotimes (tK_A) \hotimes t{\rm det}\sN^*_{A \vert X}) =
0$$
Thus the claim results from our basic inequality. \qed

\bigskip \noindent Notice that in (4.7) we do not claim that $X$ is
uniruled. If $X$ and $A$ are both smooth, this can be proved:

\begin{theorem} Assume in (4.7) additionally that $X$ and $A$ are smooth.
Then $X$ is uniruled.
\end{theorem}

\proof Let $(C_t)$ be a covering family of rational curves. Then $T_A \vert
C_t$ is nef
for general $t.$ This is well known and is easily seen by considering the
graph of the
family. Now the normal bundle sequence and our assumption on the normal
bundle imply that
$T_X \vert C_t$ is nef for general $t.$ Hence $X$ is uniruled, since the
deformations of $C_t$
fill up $X$ (if $C_t$ consider the normalization $f_t : {\Bbb P} _1 \mapsto
C_t$ and deform the
morphism $f_t,$ compare e.g., \cite{Ko}). \qed

\bigskip \noindent The main difficulty for proving (4.7) without some
hypotheses on the
singularities of $X$ and $A$ is the lacking of a dimension estimate of the
space of maps
${\rm Hom}(\Bbb P _1,X)$, cp. \cite{Ko}.

\begin{corollary} Let $(X,A)$ be a normal projective pair. Assume $\dim A
\leq 3$ and that
$A$ has only terminal singularities. Assume that $\sym m \sN_{A \vert X}
\simeq A \hotimes B$
with $A$ nef and $B$ $\Bbb Q$-effective. Then $$\grk (X) \leq \grk (A) +
\cod_X A.$$ \end{corollary}

\bigskip \noindent Of course we conjecture that all the above results hold
without
restriction on $\dim A.$ Here are some more special cases

\begin{proposition} Let $(X,A)$ be a normal projective pair. Then $\grk (X)
\leq \grk (A) +
\cod_X A,$ if one of the following conditions holds. \begin{enumerate}
\item $\sN_{A \vert X}$ is $\Bbb R$-effective and $\dim A = 1$; \item
$\sN_{A \vert X} = A \hotimes B$ with a nef $\Bbb Q$-divisor $A$ and an
effective $\Bbb Q$-divisor $B$,
moreover $X$ is $\Bbb Q$-Gorenstein of dimension $n,$ $\dim A = n-1$ and
$\grk (X) = n.$
\end{enumerate} \end{proposition}

\proof (1) First suppose $\dim A = 1.$ Then $X$ is smooth in a neighborhood
of $A,$ and the normal
bundle $\sN_A$ is nef. By adjunction we have $$ K_X \cdot A = 2g - 2 + {\rm
deg}(\sN^*_A),$$ and by the nefness of $\sN_A$ we get $K_X \cdot A \leq 2g
- 2, g = g(A).$ If $g = 0,$ we obtain $K_X \cdot A < 0$ and therefore $\grk
(X) = - \infty$ by (3.6).
If $g = 1$ we have $K_X \cdot A \leq 0 $ and we conclude by (3.7). If $g
\geq 2,$ there
is nothing to prove.

(2) By lemma (4.14) below, we have $\kappa (X \vert A) = n-1.$ From the
adjunction formula it now
follows easily that $\grk (A) = n-1.$

\bigskip \noindent The following result is easily verified by the reader
along the lines of the results of this section.

\begin{theorem} Let $(X,A)$ be a normal projective pair with $A \subset X$
Cartier. Assume that
$\sN_{A \vert X}$ is in the closure of the effective cone, that $A$ has
only terminal singularities and that
$A$ has a good minimal model via contractions and flips. Then $$\grk (X) \leq \grk (A) + 1.$$
\end{theorem}

\bigskip \noindent In the rest of this section we are discussing
applications of (4.1) and
related results.

\begin{theorem}Let $X$ be an normal projective variety. Let $\sE$ be a
vector bundle of rank $r$ on $X$. Assume that there is a section $s$ of
$\sE$ that vanishes on an normal variety $A$ such that: \begin{enumerate}
\item on $X-\sing X$, $s$ is transverse regular; and \item $\cod_A(A\cap
sing X)\ge 2$.
\end{enumerate}
If $\sE$ is ample or generically spanned then $\grk(X)\le \grk(A)+\rank
\sE$. \end{theorem}
\proof Just observe that $\sN_{A \vert X} = \sE \vert A$ is $\Bbb
Q$-effective. \qed

\bigskip \noindent Of course similar theorems can be stated with $\Bbb
R$-effectivity or nefness conditions on $\sE,$ e.g., $\sE$ is nef, $A$ has
only terminal singularities and $\dim A \leq 3.$ The last condition can be
restated as
$\rank \sE \geq \dim X - 3.$

\begin{theorem}Let $(X,A)$ be a normal projective pair. If $\grk(\det
\normB A X)>\grk(A)$ and if $\normB A X$ is $\rat$-effective, then
$\grk(X)=-\infty$. \end{theorem}
\proof If
$\grk(X)\ge 0$ then by (2.1) there would be a not identically zero section
of $kK_{X \vert A}
\hotimes \sym t\conormB A X$ for some $k>0$ and some $t\ge 0$. Since
$\normB A X$ is $\rat$-effective we conclude that we have a not identically
zero section of $kK_{X \vert A}$ for some $k>0$. This is absurd because the
condition $\grk(\det \normB A X)>\grk(A)$ implies that $\grk(K_{X \vert
A})=-\infty$ by adjunction. \qed

As an application consider the situation when we have a connected complex
submanifold $A$ of a projective manifold $X$ with $\normB A X$ ample. If
$\grk(A)\not=\dim A$ then we conclude that $\grk(X)=-\infty$. The ampleness
of
$\normB A X$ can be reduced in certain situations.

\begin{lemma}Let $(X,A)$ be a normal projective pair. Assume that $A\subset
X$ is a divisor whose normal sheaf can be written as $\Bbb Q$-Weil divisor
as a sum of a nef and an effective divisor. If $\sL$ is a rank one
reflexive free sheaf on $X$ adapted to $(X,A)$ with $\grk(\sL) =\dim X$,
then $\grk(\sL_A) =\dim A$. \end{lemma}
\proof
Write (as $\Bbb Q-$divisors) $L = H + E$ with $H$ ample and $E$ effective.
Now add all components $E_i$ of
$E$ to $H$ which have the property that $E_i \vert A$ is the sum a nef and
an effective $\Bbb Q$-divisor.
Then we obtain a decomposition
$$ L = H' + E'$$
and no component of $E'$ has this property. Therefore $A$ is not a
component of $E'.$ Thus $H' \vert A$ is
big and $E' \vert A$ is effective, so that $L \vert A $ is big.

\begin{corollary}
Let $(X,A)$ be a normal pair with $A\subset X$ be a divisor whose normal
sheaf is
(as $\Bbb Q$-Weil divisor) a sum of an effective and a nef divisor. If $X$
is of general
type, then A is of general type, also.
\end{corollary}

\bigskip \noindent Our last application concerns covering families of
subvarieties.

\begin{theorem} Let $(X,A)$ be a normal projective pair. Assume that $tA$
moves as a cycle
in a covering family in $X$ for some $t > 0.$ Then $$ \grk (X) \leq \grk
(A) + \cod_X A.$$
\end{theorem}

\proof We want to prove that $\sN_{A \vert X}$ is $\rat$-effective. Going
to the graph of
the covering family, we see immediately that $\sN_{tA \vert X}$ is
generically spanned.
Via the canonical, generically injective map $\sI^t_A/\sI^{t+1}_A \to
S^t(\sI/\sI^2)$ it
follows that $\sym t \sN_A$ is generically spanned, hence $\sN_A$ is $\Bbb
Q$-effective.

\begin{corollary} Let $X$ be a projective manifold, $(A_t)$ a covering
family of positive
dimensional subvarieties. Assume that $A_t$ is smooth for general $t.$ If
$A_t$ is not
of general type ($t$ general), then $X$ is not of general type. \end{corollary}

\vspace{1cm}
\small
\begin{tabular}{lcl}
Thomas Peternell&& Andrew J. Sommese \\ Mathematisches Institut & &
Department of Mathematics\\ Universit\" at Bayreuth & & University of Notre
Dame\\ D-95440 Bayreuth, Germany&&Notre Dame, Indiana 46556, U.S,A,\\
thomas.peternell@uni-bayreuth.de && sommese.1@nd.edu\\ &&URL: {\tt
http://www.nd.edu/$\sim$sommese/index.html}\\ fax: Germany + 921--552999&&
fax: U.S.A. + 219--631-6579 \\ \end{tabular}

\end{document}